\documentclass{amsart}
\usepackage{amssymb}
\theoremstyle{plain}
 \newtheorem{thm}{Theorem}
 \newtheorem{cor}{Corollary}
 \newtheorem{lem}{Lemma}
 \newtheorem{prop}{Proposition}
\theoremstyle{definition}

\theoremstyle{remark}
 \newtheorem{rem}{Remark\ignorespaces}

\allowdisplaybreaks
\newcommand{\NaturalNumber}{\mathbb N}
\newcommand{\RealNumber}{\mathbb R}
\renewcommand{\labelenumi}{(\roman{enumi})}

\begin{document}
\title[Commutative Semigroups of Nonexpansive Mappings]
{Common Fixed Points of Commutative Semigroups of Nonexpansive Mappings}
\author[T. Suzuki]{Tomonari Suzuki}
\date{}
\hyphenation{kyu-shu kita-kyu-shu to-bata-ku sen-sui-cho}
\address{
Department of Mathematics,
Kyushu Institute of Technology,
1-1, Sensuicho, Tobataku, Kitakyushu 804-8550, Japan}
\email{suzuki-t@mns.kyutech.ac.jp}
\keywords{Nonexpansive semigroup, Common fixed point,
 Invariant mean, Opial property}
\subjclass[2000]{Primary 47H20, Secondary 47H10}

\begin{abstract}
In this paper,
 we discuss characterizations of common fixed points
 of commutative semigroups of nonexpansive mappings.
We next prove convergence theorems to a common fixed point.
We finally discuss
 nonexpansive retractions onto the set of common fixed points.
In our discussion,
 we may not assume the strict convexity of the Banach space.
\end{abstract}
\maketitle

\section{Introduction}
\label{SC:introduction}

Let $C$ be a closed convex subset of a Banach space $E$.
A mapping $T$ on $C$ is called a {\it nonexpansive mapping} if
 $ \| Tx - Ty \| \leq \| x - y \| $
 for all $x,y \in C$.
We denote by $F(T)$ the set of fixed points of $T$.
Kirk \cite{REF:Kirk1965_AMMonth} proved that
 $F(T)$ is nonempty
 in the case that $C$ is weakly compact and has normal structure.
See also
 \cite{REF:Baillon1979_1,
 REF:Browder1965_ProcNAS_4,
 REF:Browder1965_ProcNAS_3,
 REF:Gohde1965}
 and others.
If $C$ is weakly compact and $E$ has the Opial property,
 then $C$ has normal structure; see \cite{REF:Gossez_LamiDazo1972_Pacific}.
Thus, $F(T)$ is nonempty
 in the case that $C$ is weakly compact and $E$ has the Opial property.

Let $(S,+)$ be a {\it commutative semigroup}, i.e.,
\begin{enumerate}
\item $s + t \in S$ for $s, t \in S$;
\item $(s + t) + u = s + (t + u)$ for $s, t, u \in S$; and
\item $s + t = t + s$ for $s, t \in S$.
\end{enumerate}
Then a family $\{ T(t) : t \in S \}$ of mappings on $C$ is called
 a {\it commutative semigroup of nonexpansive mappings} on $C$
 (a {\it nonexpansive semigroup} on $C$, for short)
 if the following are satisfied:
 \begin{enumerate}
 \renewcommand{\labelenumi}{(sg \arabic{enumi})}
 \item
 for each $t \in S$,
 $T(t)$ is a nonexpansive mapping on $C$; and
 \item
 $ T(s+t) = T(s) \circ T(t) $ for all $s, t \in S$.
 \end{enumerate}
We put $F({\mathcal S}) = \bigcap_{t \in S} F \big( T(t) \big)$.
Common fixed point theorems for families of nonexpansive mappings
 are proved in
 \cite{REF:Browder1965_ProcNAS_3,
 REF:Bruck1974_Pacific,
 REF:DeMarr1963_Pacific} and others.
The following is the corollary of the famous theorem proved
 by Bruck \cite{REF:Bruck1974_Pacific}.

\begin{thm}[Bruck \cite{REF:Bruck1974_Pacific}]
\label{THM:Bruck}
Let $S$ be a commutative semigroup and
 let $\{ T(t) : t \in S \}$ be a nonexpansive semigroup
 on a weakly compact convex subset $C$ of a Banach space $E$.
Suppose that $C$ has the fixed point property for nonexpansive mappings.
Then $F({\mathcal S})$ is a nonempty
 nonexpansive retract of $C$.
\end{thm}

\noindent
We note that from this theorem,
 $F({\mathcal S})$ is nonempty
 in the case that $C$ is weakly compact and $E$ has the Opial property.

Many convergence theorems for nonexpansive mappings and
 families of nonexpansive mappings have been studied;
 see \cite{REF:Atsushiba_Shioji_Takahashi2000_JNCA_2,
 REF:Baillon1975,
 REF:Browder1967_ARMA,
 REF:Edelstein1966_AMMonth,
 REF:Halpern1967_BullAMS,
 REF:KrasnoselskiiMA1955_Uspekhi,
 REF:Linhart1972_MonatMath,
 REF:Reich1979_JMAA,
 REF:TS2003_Nihonkai,
 REF:Wittmann1992_ArchMath} and others.
In these theorems, we assume the strict convexity of the Banach space $E$.
In the results of
 \cite{REF:Edelstein_Obrien1978_JLondon,
 REF:Ishikawa1976_ProcAMS,
 REF:Ishikawa1979_Pacific,
 REF:TSP_ishiinfi_2_04},
 we may not assume the strict convexity of $E$.
Very recently,
 the author in \cite{REF:TS2002_JNCA} proved
 strong convergence theorems
 for families of nonexpansive mappings
 without the assumption of the strict convexity.
See also \cite{REF:TSP_mcbt_1_05}.

In this paper,
 we extend the results in \cite{REF:TS2002_JNCA,REF:TSP_mcbt_1_05}
 to commutative semigroups of nonexpansive mappings.
We also discuss characterizations of common fixed points
 of nonexpansive semigroups
 and
 nonexpansive retractions onto the set of common fixed points.
In our discussion,
 we may not assume the strict convexity of the Banach space.

\section{Preliminaries}
\label{SC:preliminaries}

Throughout this paper
 we denote by $\NaturalNumber$ the set of all positive integers
 and by $\RealNumber$ the set of all real numbers.
Let $A$ be a subset of a set $S$.
Then we define a function $I_A$ from $S$ into $\RealNumber$ by
 $$ I_A (t) =
 \begin{cases}
  1, & \text{if $t \in A$}, \\
  0, & \text{if $t \not\in A$}.
 \end{cases} $$

Let $E$ be a Banach space.
We denote by $E^\ast$ the dual of $E$.
We recall that $E$ is said to have the {\it Opial property}
 \cite{REF:Opial1967_BullAMS}
 if for each weakly convergent sequence $\{ x_n \}$ in $E$
 with weak limit $x_0$,
 $\liminf_n \| x_n - x_0 \| < \liminf_n \| x_n - x \|$
 for all $x \in E$ with $x \neq x_0$.
All Hilbert spaces, all finite dimensional Banach spaces
 and $ \ell^p ( 1 \leq p < \infty ) $
 have the Opial property.
A Banach space with a duality mapping
 which is weakly sequentially continuous
 also has the Opial property; see \cite{REF:Gossez_LamiDazo1972_Pacific}.
We know that every separable Banach space can be equivalently renormed
 so that it has the Opial property; see \cite{REF:vanDulst1982_JLondon}.
See also
 \cite{REF:LamiDazo1973_ProcAMS,
 REF:LinPK_Tan_Xu1995_NATMA,
 REF:PrusS1992_NATMA,
 REF:Sims1985_1} and others.

We recall that a closed convex subset $C$ of a Banach space $E$ has
 the fixed point property for nonexpansive mappings
 if
 for every bounded closed convex subset $D$ of $C$ and
 for every nonexpansive mapping on $D$ has a fixed point.
Every compact convex subset of any Banach space
 has the fixed point property for nonexpansive mappings.
Also, does every weakly compact convex subset of
 a Banach space with the Opial property.

Let $(S,+)$ be a commutative semigroup.
Then we can consider that
 $S$ is a directed set with relation $\leq$ defined as follows:
$s \leq t$ if and only if
 $s = t$ or there exists $u \in S$ such that $s + u = t$.
We denote by $B(S)$
 the Banach space consisting of
 all bounded functions from $S$ into $\RealNumber$
 with supremum norm.
For $s \in S$, we define a mapping $\ell_s$ on $B(S)$ by
 $$ (\ell_s a) (t) = a(s+t) $$
 for $a \in B(S)$ and $t \in S$.
Let $X$ be a linear subspace of $B(S)$
 such that $I_S \in X$ and $X$ is $\ell_s$-invariant for all $s \in S$.
We call that $\mu \in X^\ast$ is a {\it mean} on $X$
 if $\| \mu \| = \mu(I_S) = 1$.
We know that $\mu$ is a mean on $X$ if and only if
 $$ \inf_{t \in S} a(t) \leq \mu(a) \leq \sup_{t \in S} a(t) $$
 for all $a \in X$;
 see \cite{REF:Takahashi_ybook} and others.
We also know that
 if $a, b \in X$ satisfies $a(t) \leq b(t)$ for all $t \in S$,
 then $\mu(a) \leq \mu(b)$.
Sometimes, we denote by $\mu_t \big( a(t) \big)$ the value $\mu(a)$.
A mean $\mu$ on $X$ is called {\it invariant}
 if
 $$ \mu_t \big( a(t) \big) = \mu_t \big( a(s+t) \big) $$
 for all $a \in X$ and $s \in S$.
We note that since $S$ is commutative,
 there exists an invariant mean on $X$.
Let $\{ \mu_\alpha : \alpha \in D \}$ be a net of means on $X$.
Then $\{ \mu_\alpha \}$ is called
 {\it asymptotically invariant} \cite{REF:Day1959_Illinois}
 if
 $$ \lim_{\alpha \in D} \Big(
 \mu_\alpha (a) - (\mu_\alpha)_t \big( a(s+t) \big) \Big) = 0 $$
 for all $a \in X$ and $s \in S$.
It is obvious that
 if $\{ \mu_\alpha \}$ is an asymptotically invariant net
 of means on $X$,
 then so is every subnet $\{ \mu_{\alpha_\beta} \}$ of $\{ \mu_\alpha \}$.

Let $E$ be a Banach space and
 let $C$ be a weakly compact convex subset of $E$.
Let $S$ be a commutative semigroup and
 let $\{ T(t) : t \in S \}$ be a nonexpansive semigroup on $C$.
Let $X$ be a linear subspace of $B(S)$
 such that $I_S \in X$, $X$ is $\ell_s$-invariant for all $s \in S$, and
 $\Big( t \mapsto f \big( T(t)x \big) \Big) \in X$
 for all $x \in C$ and $f \in E^\ast$.
Let $\mu$ be a mean on $X$.
Then we know that
 for each $x \in C$, there exists a unique element $x_0$ of $C$ satisfying
 $$ \mu_t \Big( f \big( T(t)x \big) \Big) = f(x_0) $$
 for all $f \in E^\ast$;
 see \cite{REF:Hirano_Kido_Takahashi1988_NATMA,REF:Takahashi1981_ProcAMS}.
Following Rod\'e \cite{REF:Rode1982_JMAA},
 we denote such $x_0$ by $T_\mu x$.
We also know that $T_\mu$ is a nonexpansive mapping on $C$.

We now prove the following,
 which are used in Section \ref{SC:characterization}.

\begin{lem}
\label{LEM:liminf}
Let $S$ be a commutative semigroup, and
 let $\{ \alpha_t : t \in S \}$ be a real net.
Then
 $$ \liminf_{t \in S} \alpha_t = \liminf_{t \in S} \alpha_{s+t}
 \quad\text{and}\quad
 \limsup_{t \in S} \alpha_t = \limsup_{t \in S} \alpha_{s+t}
 $$
 for $s \in S$.
\end{lem}

\begin{proof}
Fix $\lambda \in \RealNumber$ with $\liminf_{t \in S} \alpha_t < \lambda$.
For $t_1 \in S$,
 since $s + t_1 \in S$,
 there exists $t_2 \in S$ such that $t_2 \geq s + t_1$ and
 $\alpha_{t_2} < \lambda$.
In the case of $t_2 = s + t_1$,
 we put $t_3 = t_1$.
In the case that there exists $t_4 \in S$ such that
 $t_2 = s + t_1 + t_4$,
 we put $t_3 = t_1 + t_4$.
In both cases,
 we have
 $$ t_3 \geq t_1
 \quad\text{and}\quad
 \alpha_{s + t_3} = \alpha_{t_2} < \lambda . $$
Therefore
 $\liminf_{t \in S} \alpha_{s+t} \leq \lambda$.
Since $\lambda \in \RealNumber$ is arbitrary,
 we obtain
 $$ \liminf_{t \in S} \alpha_t \geq \liminf_{t \in S} \alpha_{s+t} . $$
Fix $\lambda \in \RealNumber$ with $\liminf_{t \in S} \alpha_t > \lambda$.
Then there exists $t_5 \in S$ such that
 $\alpha_t > \lambda$ for $t \geq t_5$.
Since $t \geq t_5$ implies $s + t \geq t_5$,
 we have $\alpha_{s+t} > \lambda$ for $t \geq t_5$.
Therefore
 $\liminf_{t \in S} \alpha_{s+t} \geq \lambda$.
Since $\lambda \in \RealNumber$ is arbitrary,
 we obtain
 $$ \liminf_{t \in S} \alpha_t \leq \liminf_{t \in S} \alpha_{s+t} . $$
Hence
 $ \liminf_{t \in S} \alpha_t = \liminf_{t \in S} \alpha_{s+t} $.
We also have
 $$ \limsup_{t \in S} \alpha_t
 = - \liminf_{t \in S} (- \alpha_t)
 = - \liminf_{t \in S} (- \alpha_{s+t})
 = \limsup_{t \in S} \alpha_{s+t} .$$
This completes the proof.
\end{proof}

\begin{lem}
\label{LEM:S-set}
Let $S$ be a commutative semigroup, and
 let $\tilde\mu$ be a mean on $B(S)$.
Let $A_1, A_2, A_3, \cdots, A_k$ be subsets of $S$.
Put
 $$ A = \bigcap_{j=1}^k A_j
 \quad\text{and}\quad
 \alpha = \sum_{j=1}^{k} \liminf_{s \in S}
  \tilde\mu_t \big( I_{A_j}(s+t) \big) - k + 1 .$$
Suppose $\alpha > 0$.
Then
 $$ \liminf_{s \in S} \tilde\mu_t \big( I_A(s+t) \big) \geq \alpha
 \quad\text{and}\quad
 \{ s_0 + t : t \in S \} \cap A \neq \varnothing
 $$
 hold for all $s_0 \in S$.
\end{lem}

\begin{proof}
It is obvious that
 $t \in A$
 if and only if
 $ \sum_{j=1}^{k} I_{A_j} (t) = k $,
 and
 $t \in S \setminus A$
 if and only if
 $ \sum_{j=1}^{k} I_{A_j} (t) \leq k-1 $.
Therefore we obtain
 $$ I_A(t) \geq \sum_{j=1}^{k} I_{A_j} (t) - k + 1 $$
 for all $t \in S$.
Hence,
 \begin{align*}
 \liminf_{s \in S} \tilde\mu_t \big( I_A(s+t) \big)
 &\geq \liminf_{s \in S} \tilde\mu_t
   \left( \sum_{j=1}^{k} I_{A_j} (s+t) - k + 1 \right) \\*
 &= \liminf_{s \in S} \left( \sum_{j=1}^{k} \tilde\mu_t
   \big( I_{A_j} (s+t) \big) - k + 1 \right) \\
 &\geq \sum_{j=1}^{k} \liminf_{s \in S} \tilde\mu_t
   \big( I_{A_j} (s+t) \big) - k + 1  \\*
 &= \alpha > 0 .
 \end{align*}
Therefore there exists $s_1 \in S$ such that
 $$ \inf_{s \geq s_1} \tilde\mu_t \big( I_A(s+t) \big)
 \geq \frac{\alpha}{2} . $$
We suppose that there exists $s_0 \in S$ such that
 $\{ s_0 + t : t \in S \} \cap A = \varnothing$.
Then
 $\{ s_0 + s_1 + t : t \in S \} \cap A = \varnothing$
 and hence
 $I_A(s_0+s_1+t) = 0$ for $t \in S$.
Since $s_0 + s_1 \geq s_1$,
 we obtain
 $$ 0
 < \frac{\alpha}{2}
 \leq \tilde\mu_t \big( I_A(s_0+s_1+t) \big)
 = \tilde\mu(0)
 = 0 . $$
This is a contradiction.
This completes the proof.
\end{proof}

\begin{lem}
\label{LEM:S-set-invmean}
Let $S$ be a commutative semigroup and
 let $\mu$ be an invariant mean on $B(S)$.
Let $A_1, A_2, A_3, \cdots, A_k$ be subsets of $S$.
Put
 $$ A = \bigcap_{j=1}^k A_j
 \quad\text{and}\quad
 \alpha = \sum_{j=1}^{k} \mu(I_{A_j}) - k + 1 .$$
Suppose $\alpha > 0$.
Then
 $ \mu(I_A) \geq \alpha $
 holds and
 $ \{ s_0 + t : t \in S \} \cap A \neq \varnothing $
 hold for all $s_0 \in S$.
\end{lem}

\begin{proof}
For every subset $B$ of $S$, we have
 $$ \liminf_{s \in S} \mu_t \big( I_B(s+t) \big)
 = \liminf_{s \in S} \mu_t \big( I_B(t) \big)
 = \mu ( I_B ) . $$
So, by Lemma \ref{LEM:S-set},
 we obtain the desired result.
\end{proof}

\section{Characterizations}
\label{SC:characterization}

In this section,
 we discuss the characterization of common fixed points.

We first prove the following,
 which plays an important role in this paper.

\begin{thm}
\label{THM:seq}
Let $E$ be a Banach space and
 let $C$ be a weakly compact convex subset of $E$.
Let $S$ be a commutative semigroup and
 let $\{ T(t) : t \in S \}$ be a nonexpansive semigroup on $C$.
Let $X$ be a linear subspace of $B(S)$
 such that $I_S \in X$, $X$ is $\ell_s$-invariant for all $s \in S$, and
 $\Big( t \mapsto f \big( T(t)x \big) \Big) \in X$
 for all $x \in C$ and $f \in E^\ast$.
Let $\mu$ be an invariant mean on $X$.
Suppose that $T_\mu z = z$ for some $z \in C$.
Then there exist sequences $\{ p_n \}$ and $\{ q_n \}$ in $S$
 and $\{ f_n \}$ in $E^\ast$ such that
 \begin{align*}
 & p_{n+1} = p_n + q_n, \\*
 & \| T(p_n) z - z \| \geq \lambda - \frac{1}{3^{n+1}}, \\
 & f_\ell \big( T(p_n) z - z \big) \leq \frac{2^{\ell+1}}{3^{\ell+1}}
  \quad\text{for}\; \ell = 1, 2, \cdots, n-1, \\*
 & \| f_n \| = 1 \quad\text{and}\quad
 f_n \big( T(p_n) z - z \big) = \| T(p_n) z - z \|
 \end{align*}
 for all $n \in \NaturalNumber$,
 where
 $$ \lambda = \limsup_{t \in S} \| T(t) z - z \| .$$
\end{thm}

Before proving it, we need some preliminaries.
In the following lemmas and the proof of Theorem \ref{THM:seq},
 we put
 $$ A(f,\varepsilon) =
  \{ t \in S : f \big( T(t) z - z \big) \leq \varepsilon \}$$
 for $f \in E^\ast$ and $\varepsilon > 0$,
 and
 $$ B(\varepsilon) = \{ t \in S :
  \| T(t) z - z \| \geq \lambda - \varepsilon \} $$
 for $\varepsilon > 0$.
By the Hahn-Banach theorem,
 there exists an extension $\tilde\mu$ of $\mu$ such that
 the domain of $\tilde\mu$ is $B(S)$ and $\| \tilde\mu \| = \| \mu \| = 1$.
It is obvious that such $\tilde\mu$ is a mean on $B(S)$.

\begin{lem}
\label{LEM:leq-lambda}
For every $s \in S$,
 $ \| T(s) z - z \| \leq \lambda $
 holds.
\end{lem}

\begin{proof}
Fix $s \in S$ and $\varepsilon > 0$.
Then by the definition of $\lambda$,
 there exists $t_0 \in S$ such that
 $$ \sup_{t \geq t_0} \| T(t) z - z \| \leq \lambda + \varepsilon . $$
Hence, for each $t \in S$, we have
 $$ \| T(t_0+t) z - z \| \leq \lambda + \varepsilon $$
 because $t_0+t \geq t_0$.
By the Hahn-Banach theorem,
 there exists $f \in E^\ast$ such that
 $$ \| f \| = 1
 \quad\text{and}\quad
 f \big( T(s) z - z \big) = \| T(s) z - z \| . $$
For $t \in S$,
 we have
 \begin{align*}
 \| T(s) z - z \|
 &= f \big( T(s) z - z \big) \\*
 &= f \big( T(s) z - T(s+t_0+t) z \big) + f \big( T(s+t_0+t) z - z \big) \\
 &\leq \| f \| \; \| T(s) z - T(s+t_0+t) z \|
  + f \big( T(s+t_0+t) z - z \big) \\
 &= \| T(s) z - T(s) \circ T(t_0+t) z \|
 + f \big( T(s+t_0+t) z - z \big) \\
 &\leq \| T(t_0+t) z - z \| + f \big( T(s+t_0+t) z - z \big) \\*
 &\leq \lambda + \varepsilon + f \big( T(s+t_0+t) z - z \big) .
 \end{align*}
Since $\mu$ is an invariant mean on $X$,
 we have
 \begin{align*}
 \| T(s) z - z \|
 &= \mu_t \big( \| T(s) z - z \| \big) \\*
 &\leq \mu_t \Big( \lambda + \varepsilon + f \big( T(s+t_0+t) z - z \big)
   \Big) \\
 &= \lambda + \varepsilon
  + \mu_t \Big( f \big( T(s+t_0+t) z \big) \Big) - f(z) \\
 &= \lambda + \varepsilon
  + \mu_t \Big( f \big( T(t) z \big) \Big) - f(z) \\
 &= \lambda + \varepsilon + f(T_\mu z) - f(z) \\*
 &= \lambda + \varepsilon .
 \end{align*}
Since $\varepsilon > 0$ is arbitrary,
 we have $ \| T(s) z - z \| \leq \lambda $.
This completes the proof.
\end{proof}

\begin{lem}
\label{LEM:Afe}
Fix $s_0 \in S$ and
 let $f \in E^\ast$ with
 $$\| f \| = 1
 \quad\text{and}\quad
 f \big( T(s_0) z - z \big) = \| T(s_0) z - z \| .$$
Let $\delta$ be a positive real number
 satisfying
 $ \| T(s_0) z - z \| \geq \lambda - \delta $.
Then
 $$ \liminf_{s \in S} \tilde\mu_t \big( I_{A(f,\varepsilon)} (s+t) \big)
 \geq \frac{\varepsilon}{\varepsilon+\delta} $$
 hold for all $\varepsilon > 0$.
\end{lem}

\begin{proof}
For $t \in S$,
 by Lemma \ref{LEM:leq-lambda},
 we have
 $$ \| T(s_0) z - T(s_0+t) z \|
 = \| T(s_0) z - T(s_0) \circ T(t) z \|
 \leq \| T(t) z - z \| \leq \lambda $$
 and hence
 \begin{align*}
 f \big( T(s_0+t) z - z \big)
 &= f \big( T(s_0+t) z - T(s_0) z \big) + f \big( T(s_0) z - z \big) \\*
 &\geq - \| f \| \; \| T(s_0+t) z - T(s_0) z \| + f \big( T(s_0) z - z \big) \\
 &= - \| T(s_0+t) z - T(s_0) z \| + \| T(s_0) z - z \| \\*
 &\geq - \lambda + \lambda - \delta
 = - \delta .
 \end{align*}
On the other hand,
 by the definition of $A(f,\varepsilon)$,
 $ f \big( T(t) z - z \big) > \varepsilon $
 for all $t \in S \setminus A(f,\varepsilon)$.
Therefore we have
 \begin{align*}
 f \big( T(s_0+t) z - z \big)
 &\geq - \delta \; I_{A(f,\varepsilon)} (s_0+t)
  + \varepsilon \; I_{S \setminus A(f,\varepsilon)} (s_0+t) \\*
 &= - \delta \; I_{A(f,\varepsilon)} (s_0+t)
  + \varepsilon \; I_{S} (s_0+t)
  - \varepsilon \; I_{A(f,\varepsilon)} (s_0+t) \\*
 &= - (\delta + \varepsilon) \; I_{A(f,\varepsilon)} (s_0+t) + \varepsilon
 \end{align*}
 for $t \in S$.
So, for $s \in S$, we have
 \begin{align*}
 0
 &= f(T_\mu z - z) \\*
 &= \mu_t \Big( f \big( T(t) z - z \big) \Big) \\
 &= \mu_t \Big( f \big( T(s_0+s+t) z - z \big) \Big) \\
 &= \tilde\mu_t \Big( f \big( T(s_0+s+t) z - z \big) \Big) \\
 &\geq \tilde\mu_t \Big(
   - (\delta + \varepsilon) \; I_{A(f,\varepsilon)} (s_0+s+t) + \varepsilon
  \Big) \\*
 &= - (\delta + \varepsilon) \;
  \tilde\mu_t \Big( I_{A(f,\varepsilon)}(s_0+s+t) \Big)
  + \varepsilon .
 \end{align*}
Hence, we obtain
 $$ \tilde\mu_t \Big( I_{A(f,\varepsilon)}(s_0+s+t) \Big)
 \geq \frac{\varepsilon}{\varepsilon+\delta} $$
 for all $s \in S$.
So, by Lemma \ref{LEM:liminf}, we have
 $$ \liminf_{s \in S} \tilde\mu_t \Big( I_{A(f,\varepsilon)}(s+t) \Big)
 = \liminf_{s \in S} \tilde\mu_t \Big( I_{A(f,\varepsilon)}(s_0+s+t) \Big)
 \geq \frac{\varepsilon}{\varepsilon+\delta} . $$
This completes the proof.
\end{proof}

\begin{lem}
\label{LEM:B}
$$ \liminf_{s \in S} \tilde\mu_t \big( I_{B(\varepsilon)}(s+t) \big) = 1 $$
 hold for all $\varepsilon > 0$.
\end{lem}

\begin{proof}
We fix $\varepsilon > 0$ and $\eta \in \RealNumber$ with $1/2 < \eta < 1$
 and put
 $ \delta = \varepsilon (1-\eta) / (2 \eta) $.
We note that $0 < \delta < \varepsilon/2$.
By the definition of $\lambda$,
 there exists $s_0 \in S$ such that
 $ \| T(s_0) z - z \| \geq \lambda - \delta $.
Fix $f \in E^\ast$ with
 $\| f \| = 1 $ and $ f \big( T(s_0) z - z \big) = \| T(s_0) z - z \| $.
So, by Lemma \ref{LEM:Afe}, we have
 $$ \liminf_{s \in S} \tilde\mu_t \big( I_{A(f,\varepsilon/2)}(s+t) \big)
 \geq \frac{\varepsilon/2}{\varepsilon/2+\delta}
 = \eta .$$
For $t \in S$
 with $s_0+t \in A(f,\varepsilon/2)$,
 we have
 \begin{align*}
 \| T(t) z - z \|
 &\geq \| T(s_0) z - T(s_0) \circ T(t) z \| \\*
 &= \| f \| \; \| T(s_0) z - T(s_0+t) z \| \\
 &\geq f \big( T(s_0) z - T(s_0+t) z \big) \\
 &= f \big( T(s_0) z - z \big) + f \big( z - T(s_0+t) z \big) \\
 &= \| T(s_0) z - z \| + f \big( z - T(s_0+t) z \big) \\
 &\geq \lambda - \delta - \frac{\varepsilon}{2} \\
 &\geq \lambda - \varepsilon
 \end{align*}
 and hence $t \in B(\varepsilon)$.
Therefore
 $ I_{B(\varepsilon)} (t) \geq I_{A(f,\varepsilon/2)} (s_0+t) $
 for all $t \in S$.
So, by Lemma \ref{LEM:liminf},
 we obtain
 \begin{align*}
 \liminf_{s \in S} \tilde\mu_t \big( I_{B(\varepsilon)}(s+t) \big)
 &\geq \liminf_{s \in S}
  \tilde\mu_t \big( I_{A(f,\varepsilon/2)} (s_0+s+t) \big) \\*
 &= \liminf_{s \in S}
  \tilde\mu_t \big( I_{A(f,\varepsilon/2)} (s+t) \big) \\*
 &\geq \eta .
 \end{align*}
Since $\eta$ is arbitrary, we obtain the desired result.
\end{proof}

\begin{proof}[Proof of Theorem \ref{THM:seq}]
By the definition of $\lambda$,
 there exists $p_1 \in S$ such that
 $ \| T(p_1) z - z \| \geq \lambda - 1/3^2 $.
Take $f_1 \in E^\ast$ with
 $\| f_1 \| = 1 $ and $ f_1 \big( T(p_1) z - z \big) = \| T(p_1) z - z \| $.
By Lemma \ref{LEM:Afe}, we have
 $$ \liminf_{s \in S} \tilde\mu_t
  \Big( I_{A \big( f_1,(2/3)^2 \big)}(s+t) \Big)
 \geq \frac{2^2}{2^2+1} .$$
We now define inductively
 sequences $\{ p_n \}$ in $S$
 and $\{ f_n \}$ in $E^\ast$.
Suppose $p_k \in S$ and $f_k \in E^\ast$ are known.
Since
 \begin{align*}
 & \liminf_{s \in S} \tilde\mu_t \big( I_{B(1/3^{k+2})} (s+t) \big)
  + \sum_{\ell=1}^k \liminf_{s \in S} \tilde\mu_t
   \big( I_{A \big(f_\ell,(2/3)^{\ell+1}\big)} (s+t) \big)
  - k \\*
 &\geq 1 + \sum_{\ell=1}^k \frac{2^{\ell+1}}{2^{\ell+1}+1} -k \\
 &\geq 1 + \sum_{\ell=1}^k \frac{2^{\ell+1}-1}{2^{\ell+1}} -k
 = 1 + \sum_{\ell=1}^k \frac{-1}{2^{\ell+1}} \\*
 &> \frac12 > 0 ,
 \end{align*}
 we have
 $$ \{ p_{k} + t : t \in S \} \cap B(1/3^{k+2}) \cap
  \bigcap_{\ell=1}^k A\big(f_\ell,(2/3)^{\ell+1}\big)
 \neq \varnothing $$
 by Lemma \ref{LEM:S-set}.
So we can choose $p_{k+1} \in S$ such that $ p_{k+1} = p_k + t $
 for some $t \in S$,
 $$ \| T(p_{k+1}) z - z \| \geq \lambda - \frac{1}{3^{k+2}} ,
 \quad\text{and}\quad
 f_\ell \big( T(p_{k+1}) z - z \big) \leq \frac{2^{\ell+1}}{3^{\ell+1}} $$
 for $\ell = 1, 2, \cdots, k$.
Take $f_{k+1} \in E^\ast$ with
 $$\| f_{k+1} \| = 1 \quad\text{and}\quad
 f_{k+1} \big( T(p_{k+1}) z - z \big) =
  \| T(p_{k+1}) z - z \| .$$
Note that
 $$ \liminf_{s \in S} \tilde\mu_t
  \big( I_{A\big(f_{k+1},(2/3)^{k+2}\big)} (s+t) \big)
 \geq \frac{2^{k+2}}{2^{k+2}+1} $$
 by Lemma \ref{LEM:Afe}.
Hence we have defined $\{ p_n \}$ and $\{ f_n \}$.
For each $n \in \NaturalNumber$,
 there exists $t \in S$ such that $p_{n+1} = p_n + t$.
We put $q_n = t$.
So we have defined a sequence $\{ q_n \}$ in $S$.
\end{proof}

Now, we prove the following characterization.

\begin{thm}
\label{THM:chara-Opial}
Let $C$ be a weakly compact convex subset of
 a Banach space $E$ with the Opial property.
Let $S$, $\{ T(t) : t \in S \}$, $X$ and $\mu$ be as in Theorem \ref{THM:seq}.
Then
 $z \in C$ is a common fixed point of $\{ T(t) : t \in S \}$
 if and only if
 $T_\mu z = z$.
\end{thm}

\begin{proof}
We assume that $z$ is a common fixed point of $\{ T(t) : t \in S \}$.
Then for $f \in E^\ast$, we have
 $$ f(T_{\mu} z)
 = \mu_t \; \Big( f \big( T(t) z \big) \Big)
 = \mu_t \; \Big( f(z) \Big)
 = f(z) .$$
Hence we obtain $T_{\mu}z = z$.
Conversely, we assume $T_\mu z = z$.
By Theorem \ref{THM:seq}, there exist
 sequences $\{ p_n \}$ and $\{ q_n \}$ in $S$
 and $\{ f_n \}$ in $E^\ast$ satisfying the conclusion
 of Theorem \ref{THM:seq}.
We put $\lambda = \limsup_{t \in S} \| T(t) z - z \|$.
Since $C$ is weakly compact,
 there exists a subsequence $\{ p_{n_k} \}$ of $\{ p_n \}$
 such that
 $\{ T(p_{n_k}) z \}$ converges weakly to some point $u \in C$.
If $n_k > \ell$, then
 $$ f_\ell \big( T(p_{n_k}) z - z \big) \leq \frac{2^{\ell+1}}{3^{\ell+1}} .$$
So we obtain
 $$ f_\ell(u - z) \leq \frac{2^{\ell+1}}{3^{\ell+1}} $$
 for all $\ell \in \NaturalNumber$.
Since
 \begin{align*}
 \| T(p_\ell) z - u \|
 &= \| f_\ell \| \; \| T(p_\ell) z - u \| \\*
 &\geq f_\ell \big( T(p_\ell) z - u \big) \\
 &= f_\ell \big( T(p_\ell) z - z \big) + f_\ell(z - u) \\
 &= \| T(p_\ell) z - z \| + f_\ell(z - u) \\*
 &\geq \lambda - \frac{1}{3^{\ell+1}} - \frac{2^{\ell+1}}{3^{\ell+1}}
 \end{align*}
 for $\ell \in \NaturalNumber$, we have
 $ \liminf_{\ell} \| T(p_\ell) z - u \|
 \geq \lambda $.
By Lemma \ref{LEM:leq-lambda}, we have
 \begin{align*}
 \liminf_{k \rightarrow \infty} \| T(p_{n_k}) z - z \|
 &\leq \lambda \\*
 &\leq \liminf_{\ell \rightarrow \infty} \| T(p_\ell) z - u \| \\*
 &\leq \liminf_{k \rightarrow \infty} \| T(p_{n_k}) z - u \| .
 \end{align*}
By the Opial property of $E$,
 we obtain $z = u$.
Using Lemma \ref{LEM:leq-lambda} again,
 for each $\ell \in \NaturalNumber$,
 we also have
 \begin{align*}
 & \liminf_{k \rightarrow \infty} \| T(p_{n_k}) z - T(p_\ell) z \| \\*
 &= \liminf_{k \rightarrow \infty}
   \| T(p_\ell) \circ
    T(q_\ell + q_{\ell+1} + q_{\ell+2} + \cdots + q_{n_k-1}) z
   - T(p_\ell)z \| \\
 &\leq \liminf_{k \rightarrow \infty}
   \| T(q_\ell + q_{\ell+1} + q_{\ell+2} + \cdots + q_{n_k-1}) z - z \| \\*
 &\leq \lambda .
 \end{align*}
By the Opial property of $E$,
 we obtain $T(p_\ell)z = u = z$ for all $\ell \in \NaturalNumber$.
Therefore $\lambda = 0$.
By Lemma \ref{LEM:leq-lambda}, we obtain $T(t)z = z$ for all $t \in S$.
This completes the proof.
\end{proof}

Using Theorem \ref{THM:chara-Opial}, we obtain another characterization.

\begin{thm}
\label{THM:chara2-Opial}
Let $C$ be a weakly compact convex subset of
 a Banach space $E$ with the Opial property.
Let $S$, $\{ T(t) : t \in S \}$ and $X$ be as in Theorem \ref{THM:seq}.
Let $\{ \mu_\alpha : \alpha \in D \}$ be
 an asymptotically invariant net of means on $X$.
Then
 $z \in C$ is a common fixed point of $\{ T(t) : t \in S \}$
 if and only if
 $\{ T_{\mu_\alpha} z \}$ converges weakly to $z$.
\end{thm}

\begin{proof}
We assume that $z$ is a common fixed point of $\{ T(t) : t \in S \}$.
As in the proof of Theorem \ref{THM:chara-Opial},
 we obtain $T_{\mu_\alpha}z = z$ for all $\alpha \in D$.
Therefore $\{ T_{\mu_\alpha} z \}$ converges weakly to $z$.
Conversely, we assume that $\{ T_{\mu_\alpha} z \}$ converges weakly to $z$.
Using Alaoglu's theorem,
 there exists a subnet $\{ \mu_{\alpha_\beta} : \beta \in D' \}$
 such that $\{ \mu_{\alpha_\beta} : \beta \in D' \}$ converges weakly$^\ast$
 to some point $\mu \in X^\ast$.
We know that such $\mu$ is an invariant mean on $X$;
 see \cite{REF:Day1959_Illinois}.
Since $\{ T_{\mu_\alpha} z \}$ converges weakly to $z$,
 we have
 $$ \lim_{\beta \in D'} f( T_{\mu_{\alpha_\beta}} z)
 = \lim_{\alpha \in D} f( T_{\mu_{\alpha}} z)
 = f(z) $$
 for $f \in E^\ast$.
We also have
 \begin{align*}
 \lim_{\beta \in D'} f( T_{\mu_{\alpha_\beta}} z)
 &= \lim_{\beta \in D'} (\mu_{\alpha_\beta})_t \;
   \Big( f \big( T(t)z \big) \Big) \\*
 &= \mu_t \; \Big( f \big( T(t)z \big) \Big) \\*
 &= f(T_\mu z)
 \end{align*}
 for $f \in E^\ast$.
Since $f(T_\mu z) = f(z)$ for all $f \in E^\ast$,
 we obtain $T_\mu z = z$.
By Theorem \ref{THM:chara-Opial},
 $z$ is a common fixed point of $\{ T(t) : t \in S \}$.
\end{proof}

\begin{rem}
In Theorems \ref{THM:chara-Opial} and \ref{THM:chara2-Opial},
 we may replace ``$E$ has the Opial property'' with the following condition:
For each weakly convergent sequence $\{ x_n \}$ in $C$
 with weak limit $x_0 \in C$,
 $\liminf_n \| x_n - x_0 \| < \liminf_n \| x_n - x \|$
 for all $x \in E$ with $x \neq x_0$.
We remark that if $C$ is a compact subset of any Banach space,
 then the above condition is satisfied
 even in the case that $E$ does not have the Opial property;
 see the remark of Theorem 4 in \cite{REF:TS2004_BullAuMS}.
\end{rem}

{}From the above-mentioned remark,
 we obtain the following.

\begin{thm}
\label{THM:chara-compact}
Let $C$ be a compact convex subset of a Banach space $E$.
Let $S$, $\{ T(t) : t \in S \}$, $X$ and $\mu$ be as in Theorem \ref{THM:seq}.
Let $\{ \mu_n \}$ be
 an asymptotically invariant sequence of means on $X$.
Let $z \in C$.
Then the following are equivalent:
\begin{enumerate}
\item $z$ is a common fixed point of $\{ T(t) : t \in S \}$;
\item $T_\mu z = z$;
\item $\{ T_{\mu_n} z \}$ converges strongly to $z$.
\end{enumerate}
\end{thm}

\section{Convergence Theorem}
\label{SC:convergence}

In \cite{REF:Atsushiba_Shioji_Takahashi2000_JNCA_2},
 Atsushiba, Shioji and Takahashi proved convergence theorems
 with the assumption that the Banach space is strictly convex.
In this section, we prove a convergence theorem
 to a common fixed point
 without the assumption of the strict convexity.

The following lemma concerns Krasnoselskii and Mann's type sequences
 \cite{REF:KrasnoselskiiMA1955_Uspekhi,REF:Mann1953_ProcAMS},
 and is proved in \cite{REF:TS2002_JNCA}.
See also \cite{REF:Goebel_Kirk1983_Contemp}.

\begin{lem}[\cite{REF:TS2002_JNCA} Lemma 2]
\label{LEM:bounded-seq}
Let $\{ z_n \}$ and $\{ w_n \}$ be bounded sequences in a Banach space $E$ and
 let $\{ \alpha_n \}$ be a sequence in $[0,1]$
 with $0 < \liminf_n \alpha_n \leq \limsup_n \alpha_n < 1$.
Suppose that
 $z_{n+1} = \alpha_n w_n + (1 - \alpha_n) z_n$
 for all $n \in \NaturalNumber$ and
 $$ \limsup_{n \rightarrow \infty}
 \Big( \| w_n - w_{n+k} \| - \| z_n - z_{n+k} \| \Big) \leq 0$$
 for all $k \in \NaturalNumber$.
Then
 $\liminf_{n} \| w_n - z_n \| = 0$.
\end{lem}

Using Lemma \ref{LEM:bounded-seq},
 we obtain the following convergence theorem.

\begin{thm}
\label{THM:converge-main}
Let $C$ be a compact convex subset of a Banach space $E$.
Let $S$, $\{ T(t) : t \in S \}$ and $X$ be as in Theorem \ref{THM:seq}.
Let $\{ \mu_n \}$ be
 an asymptotically invariant sequence of means on $X$.
Suppose that
 $\lim_n \| \mu_n - \mu_{n+1} \| = 0$.
Let $x_1 \in C$ and
 define a sequence $\{ x_n \}$ in $C$ by
 $$ x_{n+1} =  \alpha_n T_{\mu_n} x_n + (1 - \alpha_n) x_n $$
 for $n \in \NaturalNumber$,
 where $\{ \alpha_n \}$ is a sequence in $[0,1]$ such that
 $0 < \liminf_{n} \alpha_n \leq \limsup_{n} \alpha_n \allowbreak < 1$.
Then $\{ x_n \}$ converges strongly to
 a common fixed point $z_0$ of $\{ T(t) : t \in S \}$.
\end{thm}

\begin{proof}
For $n, k \in \NaturalNumber$,
 taking $f \in E^\ast$ with
 $$ \| f \| = 1
 \quad\text{and}\quad
 f \big( T_{\mu_n} x_{n+k} - T_{\mu_{n+k}} x_{n+k} \big)
 = \| T_{\mu_n} x_{n+k} - T_{\mu_{n+k}} x_{n+k} \| , $$
 we have
 \begin{align*}
 & \| T_{\mu_n} x_n - T_{\mu_{n+k}} x_{n+k} \| - \| x_n - x_{n+k} \| \\*
 &\leq \| T_{\mu_n} x_n - T_{\mu_n} x_{n+k} \|
  + \| T_{\mu_n} x_{n+k} - T_{\mu_{n+k}} x_{n+k} \|
  - \| x_n - x_{n+k} \| \\
 &\leq \| x_n - x_{n+k} \|
  + \| T_{\mu_n} x_{n+k} - T_{\mu_{n+k}} x_{n+k} \|
  - \| x_n - x_{n+k} \| \\
 &= \| T_{\mu_n} x_{n+k} - T_{\mu_{n+k}} x_{n+k} \| \\
 &= f \big( T_{\mu_n} x_{n+k} - T_{\mu_{n+k}} x_{n+k} \big) \\
 &= (\mu_n)_t \; \Big( f \big( T(t) x_{n+k} \big) \Big)
  - (\mu_{n+k})_t \; \Big( f \big( T(t) x_{n+k} \big) \Big) \\
 &= (\mu_n - \mu_{n+k})_t \; \Big( f \big( T(t) x_{n+k} \big) \Big) \\
 &\leq \| \mu_n - \mu_{n+k} \| \;
   \sup_{t \in S} \Big| f \big( T(t) x_{n+k} \big) \Big| \\
 &\leq \| \mu_n - \mu_{n+k} \| \;
   \sup_{t \in S} \| f \| \; \big\| T(t) x_{n+k} \big\| \\
 &\leq \| \mu_n - \mu_{n+k} \| \;
   \sup_{x \in C} \| x \| \\*
 &\leq \sum_{j=0}^{k-1} \| \mu_{n+j} - \mu_{n+j+1} \| \;
   \sup_{x \in C} \| x \| .
 \end{align*}
Hence we obtain
 $$ \limsup_{n \rightarrow \infty} \Big(
 \| T_{\mu_n} x_n - T_{\mu_{n+k}} x_{n+k} \| - \| x_n - x_{n+k} \|
 \Big) \leq 0 $$
 for all $k \in \NaturalNumber$.
By Lemma \ref{LEM:bounded-seq},
 we obtain $\liminf_n \| T_{\mu_n} x_n - x_n \| = 0$.
Since $C$ is compact, there exists a subsequence $\{ x_{n_k} \}$
 of $\{ x_n \}$ such that
 $$ \lim_{k \rightarrow \infty} \| T_{\mu_{n_k}} x_{n_k} - x_{n_k} \| = 0 $$
 and
 $\{ x_{n_k} \}$ converges strongly to some point $z_0$.
We have
 \begin{align*}
 & \limsup_{k \rightarrow \infty} \| T_{\mu_{n_k}} z_0 - z_0 \| \\*
 &\leq \limsup_{k \rightarrow \infty}
  \Big( \| T_{\mu_{n_k}} z_0 - T_{\mu_{n_k}} x_{n_k} \|
  + \| T_{\mu_{n_k}} x_{n_k} - x_{n_k} \|
  + \| x_{n_k} - z_0 \| \Big) \\*
 &\leq \limsup_{k \rightarrow \infty}
  \Big( 2 \; \| z_0 - x_{n_k} \|
  + \| T_{\mu_{n_k}} x_{n_k} - x_{n_k} \| \Big) = 0.
 \end{align*}
Since $\{ \mu_{n_k} \}$ is
 also an asymptotically invariant sequence of means on $X$,
 by Theorem \ref{THM:chara-compact},
 we have $z_0$ is a common fixed point of $\{ T(t) : t \in S \}$.
Since
 \begin{align*}
 \| x_{n+1} - z_0 \|
 &\leq \alpha_n \| T_{\mu_n} x_n - z_0 \| + (1 - \alpha_n) \| x_n - z_0 \| \\*
 &= \alpha_n \| T_{\mu_n} x_n - T_{\mu_n} z_0 \|
  + (1 - \alpha_n) \| x_n - z_0 \| \\
 &\leq \alpha_n \| x_n - z_0 \|
  + (1 - \alpha_n) \| x_n - z_0 \| \\*
 &= \| x_n - z_0 \|,
 \end{align*}
 we obtain $\lim_{n} \| x_n - z_0 \| = \lim_k \| x_{n_k} - z_0 \| = 0$.
This completes the proof.
\end{proof}

\begin{rem}
The following lemma (Lemma \ref{LEM:bounded-seq-new}) is
 a generalization of Lemma \ref{LEM:bounded-seq},
 which is useful and proved in \cite{REF:TSP_ishiinfi_2_04}.
Using Lemma \ref{LEM:bounded-seq-new},
 we can give the shorter proof of Theorem \ref{THM:converge-main}.
However, we do not use Lemma \ref{LEM:bounded-seq-new}
 because Reference \cite{REF:TSP_ishiinfi_2_04} is not yet published.
\end{rem}

\begin{lem}[\cite{REF:TSP_ishiinfi_2_04} Lemma 2]
\label{LEM:bounded-seq-new}
Let $\{ z_n \}$ and $\{ w_n \}$ be bounded sequences in a Banach space $E$ and
 let $\{ \alpha_n \}$ be a sequence in $[0,1]$
 with $0 < \liminf_n \alpha_n \leq \limsup_n \alpha_n < 1$.
Suppose that
 $z_{n+1} = \alpha_n w_n + (1 - \alpha_n) z_n$
 for all $n \in \NaturalNumber$ and
 $$ \limsup_{n \rightarrow \infty}
  \Big( \| w_{n+1} - w_n \| - \| z_{n+1} - z_n \| \Big) \leq 0 .$$
Then
 $\lim_{n} \| w_n - z_n \| = 0$.
\end{lem}

\section{Nonexpansive Retraction}
\label{SC:nonexpansive-retraction}

In this section,
 we prove the existence theorems of the nonexpansive retraction
 onto the set of common fixed points.

\begin{thm}
\label{THM:retraction}
Let $S$ be a commutative semigroup and
 let $\{ T(t) : t \in S \}$ be a nonexpansive semigroup
 on a weakly compact convex subset $C$ of a Banach space $E$.
Suppose that $C$ has the fixed point property for nonexpansive mappings.
Then
 there exists a nonexpansive retraction $Q$ from $C$ onto
 $F({\mathcal S})$
 satisfying
 $ Q \circ T(t) = T(t) \circ Q = Q$
 for all $t \in S$.
\end{thm}

\begin{proof}
By Theorem \ref{THM:Bruck}, there exists a nonexpansive retraction
 from $C$ onto $F({\mathcal S})$.
Let $\mu$ be an invariant mean on $B(S)$ and
 put $Q = P \circ T_\mu$.
Then $Q$ is a nonexpansive mapping on C,
 because so are $P$ and $T_\mu$.
For $x \in F({\mathcal S})$, we have $T_\mu x = x \in F({\mathcal S})$
 and hence $Qx = x$.
Therefore we have shown that $Q$ is a nonexpansive retraction
 from $C$ onto $F({\mathcal S})$.
So it is obvious that $T(t) \circ Q = Q$ for all $t \in S$.
Fix $s \in S$ and $x \in C$.
For $f \in E^\ast$, we have
 \begin{align*}
 f \big( T_\mu \circ T(s) x \big)
 &= \mu_t \Big( f \big( T(t) \circ T(s) x \big) \Big)
 = \mu_t \Big( f \big( T(t+s) x \big) \Big) \\*
 &= \mu_t \Big( f \big( T(t) x \big) \Big)
 = f(T_\mu x) .
 \end{align*}
Hence we have $T_\mu \circ T(s) x = T_\mu x$.
Therefore we obtain $Q \circ T(s) = Q$ for all $s \in S$.
This completes the proof.
\end{proof}

We give another nonexpansive retraction.

\begin{thm}
\label{THM:retraction-Opial}
Let $C$ be a weakly compact convex subset of
 a Banach space $E$ with the Opial property.
Let $S$, $\{ T(t) : t \in S \}$, $X$ and $\mu$ be as in Theorem \ref{THM:seq}.
Define a mapping $Q$ on $C$ as
 $$ Qx = \mathop{\text{weak-lim}}_{n \rightarrow \infty}
  \left( \frac{1}{2} T_\mu + \frac{1}{2} I \right)^n \circ T_\mu x $$
 for all $x \in C$,
 where $I$ is the identity mapping on $C$.
Then
 $Q$ is a nonexpansive retraction from $C$ onto
 $F({\mathcal S})$
 satisfying
 $ Q \circ T(t) = T(t) \circ Q = Q$
 for all $t \in S$.
Further
 if a closed convex subset $C'$ of $C$ is $T(t)$-invariant for all $t \in S$,
 then $C'$ is also $Q$-invariant.
\end{thm}

\begin{proof}
Fix $x \in C$.
Define a sequence $\{ x_n \}$ in $C$ by
 $x_1 = T_\mu x$ and
 $x_{n+1} = \frac{1}{2} T_\mu x_n + \frac{1}{2} x_n$
 for $n \in \NaturalNumber$.
Then by the result of Edelstein and O'Brien
 \cite{REF:Edelstein_Obrien1978_JLondon},
 $\{ x_n \}$ converges weakly to a fixed point of $T_\mu$.
Since $x \in C$ is arbitrary,
 $Q$ is well-defined and $Qx$ is a fixed point of $T_\mu$ for $x \in C$.
By Theorem \ref{THM:chara-Opial},
 we have $Qx \in F({\mathcal S})$ for all $x \in C$.
For $x, y \in C$, we have
\begin{align*}
 \| Q x - Q y \|
 &\leq \liminf_{n \rightarrow \infty}
  \left\| \left( \frac12 T_\mu + \frac12 I \right)^n T_\mu x -
  \left( \frac12 T_\mu + \frac12 I \right)^n T_\mu y \right\| \\*
 &\leq \liminf_{n \rightarrow \infty}
  \left\| T_\mu x - T_\mu y \right\| \\*
 &= \| T_\mu x - T_\mu y \|
 \leq \| x - y \|
\end{align*}
 and hence $Q$ is nonexpansive.
For each $x \in F({\mathcal S})$,
 we have $T_\mu x = x$ and hence $Qx = x$.
Therefore we have shown that $Q$ is a nonexpansive retraction
 from $C$ onto $F({\mathcal S})$.
So we also obtain that $T(t) \circ Q = Q$ for all $t \in S$.
As in the proof of Theorem \ref{THM:retraction},
 we have $T_\mu \circ T(t) = T_\mu$ for all $t \in S$.
So, by the definition of $Q$,
 we obtain that $Q \circ T(t) = Q $ for all $t \in S$.
We assume that
 a closed convex subset $C'$ of $C$ is $T(t)$-invariant for all $t \in S$.
Then since $C'$ is weakly compact and convex,
 we have that $C'$ is $T_\mu$-invariant.
So, by the definition of $Q$,
 $C'$ is also $Q$-invariant.
\end{proof}

{}From Ishikawa's convergence theorem \cite{REF:Ishikawa1976_ProcAMS},
 Theorem \ref{THM:chara-compact},
 and the proof of Theorem \ref{THM:retraction-Opial},
 we also obtain the following.

\begin{thm}
\label{THM:retraction-compact}
Let $C$ be a compact convex subset of a Banach space $E$.
Let $S$, $\{ T(t) : t \in S \}$, $X$ and $\mu$ be as in Theorem \ref{THM:seq}.
Define a mapping $Q$ on $C$ as in Theorem \ref{THM:retraction-Opial}.
Then the conclusion of Theorem \ref{THM:retraction-Opial} holds.
\end{thm}

\section{$S= \NaturalNumber \times \NaturalNumber$}
\label{SC:S=NN}

Using the results in Sections \ref{SC:characterization},
 \ref{SC:convergence} and \ref{SC:nonexpansive-retraction},
 we can prove many theorems.
In this section,
 we state the deduced theorems in the case of
 $S = \NaturalNumber \times \NaturalNumber$.
And in the next section,
 we state them in the case of $S = [0,\infty)$.

We first prove the following.

\begin{lem}
\label{LEM:NN}
Put $S = \NaturalNumber \times \NaturalNumber$ and
 define a sequence $\{ \mu_n \}$ of functions on $B(S)$ by
 $$ \mu_n (a) = \frac{1}{n^2} \sum_{i=1}^n \sum_{j=1}^n a(i,j) $$
 for $n \in \NaturalNumber$ and $a \in B(S)$.
Then $\{ \mu_n \}$ is an asymptotically invariant sequence
 of means on $B(S)$ and
 satisfies
 $\lim_n \| \mu_{n} - \mu_{n+1} \| = 0$.
\end{lem}

\begin{proof}
We know that $\{ \mu_n \}$ is an asymptotically invariant sequence
 of means on $B(S)$;
 see \cite{REF:Hirano_Kido_Takahashi1988_NATMA, REF:Rode1982_JMAA,
 REF:Takahashi_ybook}.
For $n \in \NaturalNumber$ and $a \in B(S)$,
 we have
 \begin{align*}
 & \big| \mu_n(a) - \mu_{n+1}(a) \big| \\*
 &= \left| \frac{1}{n^2} \sum_{i=1}^n \sum_{j=1}^n a(i,j)
  - \frac{1}{(n+1)^2} \sum_{i=1}^{n+1} \sum_{j=1}^{n+1} a(i,j) \right| \\
 &\leq \left(\frac{1}{n^2}-\frac{1}{(n+1)^2}\right)
   \sum_{i=1}^n \sum_{j=1}^n \big| a(i,j) \big| \\*
 &\quad
  + \frac{1}{(n+1)^2} \sum_{i=1}^{n} \big| a(i,n+1) \big|
  + \frac{1}{(n+1)^2} \sum_{j=1}^{n+1} \big| a(n+1,j) \big| \\
 &\leq \left(\frac{n^2}{n^2} - \frac{n^2}{(n+1)^2}
  + \frac{2 n + 1}{(n+1)^2} \right) \; \| a \|
 \end{align*}
 and hence
 $$ \limsup_{n \rightarrow \infty} \| \mu_n - \mu_{n+1} \|
 \leq \lim_{n \rightarrow \infty}
  \left(\frac{n^2}{n^2} - \frac{n^2}{(n+1)^2}
  + \frac{2 n + 1}{(n+1)^2} \right)
 = 0 . $$
This completes the proof.
\end{proof}

The following can be proved easily.

\begin{lem}
\label{LEM:NN-Rode}
Put $S = \NaturalNumber \times \NaturalNumber$ and
 define a sequence $\{ \mu_n \}$ of functions on $B(S)$
 as in Lemma \ref{LEM:NN}.
Let $C$ be a weakly compact convex subset of a Banach space $E$ and
 let $T$ and $U$ be nonexpansive mappings on $C$ with $T U = U T$.
Then
 $$ T_{\mu_n} x
 = \frac{1}{n^2} \sum_{i=1}^n \sum_{j=1}^n T^i U^j x $$
 for all $n \in \NaturalNumber$ and $x \in C$.
\end{lem}

Using Lemmas \ref{LEM:NN} and \ref{LEM:NN-Rode},
 we obtain the following.

\begin{cor}[\cite{REF:TS2002_JNCA}]
\label{COR:converge-NN}
Let $C$ be a compact convex subset of a Banach space $E$ and
 let $T$ and $U$ be nonexpansive mappings on $C$ with $T U = U T$.
Let $x_1 \in C$ and
 define a sequence $\{ x_n \}$ in $C$ by
 $$ x_{n+1} =  \frac{\alpha_n}{n^2} \sum_{i=1}^n \sum_{j=1}^n T^i U^j x_n
  +(1 - \alpha_n) x_n $$
 for $n \in \NaturalNumber$,
 where $\{ \alpha_n \}$ is a sequence in $[0,1]$ such that
 $0 < \liminf_{n} \alpha_n \leq \limsup_{n} \alpha_n \allowbreak < 1$.
Then $\{ x_n \}$ converges strongly to
 a common fixed point $z_0$ of $T$ and $U$.
\end{cor}

\begin{cor}[\cite{REF:TS2004_BullAuMS}]
\label{COR:chara-NN}
Let $E$ be a Banach space with the Opial property and
 let $C$ be a weakly compact convex subset of $E$.
Let $T$ and $U$ be nonexpansive mappings on $C$ with $T U = U T$.
Then for $z \in C$,
 the following are equivalent:
\begin{enumerate}
\item
 $z$ is a common fixed point of $T$ and~$U$;
\item
 there exists a subnet of
 a sequence
 $$ \left\{ \frac{1}{n^2} \sum_{i=1}^n \sum_{j=1}^n T^i U^j z \right\} $$
 in $C$
 converging weakly to $z$.
\end{enumerate}
\end{cor}

We also obtain the following new results.

\begin{cor}
\label{COR:chara-mu-NN-Opial}
Let $C$ be a weakly compact convex subset of
 a Banach space $E$ with the Opial property.
Let $T$ and $U$ be nonexpansive mappings on $C$ with $T U = U T$.
Let $\mu$ be an invariant mean on $B(\NaturalNumber \times \NaturalNumber)$.
Then
 $z \in C$ is a common fixed point of $T$ and $U$
 if and only if $T_\mu z = z$.
\end{cor}

\begin{cor}
\label{COR:chara-mu-NN-compact}
Let $C$ be a compact convex subset of a Banach space $E$.
Let $T$ and $U$ be nonexpansive mappings on $C$ with $T U = U T$.
Let $\mu$ be an invariant mean on $B(\NaturalNumber \times \NaturalNumber)$.
Then
 $z \in C$ is a common fixed point of $T$ and $U$
 if and only if $T_\mu z = z$.
\end{cor}

\section{$S = [0,\infty)$}
\label{SC:S=R+}

In this section,
 we state the deduced theorems in the case of
 $S = [0,\infty)$.
As in Section \ref{SC:S=NN},
 we prove the following.

\begin{lem}
\label{LEM:R+}
Put $S = [0,\infty)$ and
 let $X$ be the Banach space consisting of
 all bounded continuous functions from $S$ into $\RealNumber$
 with the supremum norm.
Let $\{ t_n \}$ be a sequence in $(0,\infty)$ satisfying
 $ \lim_{n} t_n = \infty $.
Define a sequence $\{ \mu_n \}$ of functions on $X$ by
 $$ \mu_n (a) = \frac{1}{t_n} \int_0^{t_n} a(t) \; dt $$
 for $n \in \NaturalNumber$ and $a \in X$.
Then $\{ \mu_n \}$ is an asymptotically invariant sequence
 of means on $X$.
Further, if $\{ t_n \}$ satisfies
 $\lim_{n} t_{n+1} / t_n = 1 $,
 then
 $\lim_n \| \mu_{n} - \mu_{n+1} \| = 0$.
\end{lem}

\begin{proof}
We know that $\{ \mu_n \}$ is an asymptotically invariant sequence
 of means on $X$;
 see \cite{REF:Hirano_Kido_Takahashi1988_NATMA, REF:Rode1982_JMAA,
 REF:Takahashi_ybook}.
We assume that $\lim_{n} t_{n+1} / t_n = 1 $.
For $n \in \NaturalNumber$ and $a \in X$,
 putting $ m = \min\{ t_n, t_{n+1} \} $ and
 $ M = \max\{ t_n, t_{n+1} \}$,
 we have
 \begin{align*}
 & \big| \mu_n(a) - \mu_{n+1}(a) \big| \\*
 &= \left| \frac{1}{t_n} \int_0^{t_n} a(t) \; dt
  - \frac{1}{t_{n+1}} \int_0^{t_{n+1}} a(t) \; dt \right| \\
 &= \left| \frac{1}{m} \int_0^{m} a(t) \; dt
  - \frac{1}{M} \int_0^{M} a(t) \; dt \right| \\
 &= \left| \left( \frac{1}{m} - \frac{1}{M} \right) \int_0^{m} a(t) \; dt
  - \frac{1}{M} \int_m^{M} a(t) \; dt \right| \\
 &\leq \left(\frac{1}{m}-\frac{1}{M}\right)
   \int_{0}^{m} \big| a(t) \big| \; dt
  + \frac{1}{M} \int_{m}^{M} \big| a(t) \big| \; dt \\
 &\leq \left(\frac{1}{m}-\frac{1}{M}\right)
   \int_{0}^{m} \| a \| \; dt
  + \frac{1}{M} \int_{m}^{M} \| a \| \; dt \\
 &= \left(\frac{m}{m} - \frac{m}{M}
  + \frac{M-m}{M} \right) \; \| a \| \\*
 &= \left(2 - 2 \frac{\min\{ t_n, t_{n+1} \}}{\max\{ t_n, t_{n+1} \}} \right)
  \; \| a \|
 \end{align*}
 and hence
 $$ \limsup_{n \rightarrow \infty} \| \mu_n - \mu_{n+1} \|
 \leq \lim_{n \rightarrow \infty}
  \left(2 - 2 \frac{\min\{ t_n, t_{n+1} \}}{\max\{ t_n, t_{n+1} \}} \right)
 = 0 . $$
This completes the proof.
\end{proof}

We recall that
 a family of mappings $\{ T(t) : t \geq 0 \}$ is called
 a {\it one-parameter strongly continuous semigroup of nonexpansive mappings}
 on $C$
 if the following are satisfied:
 \begin{enumerate}
 \item
 for each $t \geq 0$,
 $T(t)$ is a nonexpansive mapping on $C$;
 \item
 $ T(s+t) = T(s) \circ T(t) $ for all $s, t \geq 0$;
 \item
 for each $x \in C$,
 the mapping $t \mapsto T(t)x $ from $[0,\infty)$ into $C$ is
 strongly continuous.
 \end{enumerate}

The following can be proved easily.

\begin{lem}
\label{LEM:R+-Rode}
Let $S$, $X$, $\{ t_n \}$ and $\{ \mu_n \}$ be
 as in Lemma \ref{LEM:R+}.
Let $C$ be a weakly compact convex subset of a Banach space $E$ and
 let $\{ T(t) : t \geq 0 \}$ be
 a one-parameter strongly continuous semigroup of nonexpansive mappings on $C$.
Then
 $$ T_{\mu_n} x
 = \frac{1}{t_n} \int_{0}^{t_n} T(t) x \; dt $$
 for all $n \in \NaturalNumber$ and $x \in C$.
\end{lem}

Using Lemmas \ref{LEM:R+} and \ref{LEM:R+-Rode},
 we obtain the following.

\begin{cor}[\cite{REF:TSP_mcbt_1_05}]
\label{COR:converge-R+}
Let $C$ be a compact convex subset of a Banach space $E$ and
 let $\{ T(t) : t \geq 0 \}$ be
 a one-parameter strongly continuous semigroup of nonexpansive mappings on $C$.
Let $x_1 \in C$ and
 define a sequence $\{ x_n \}$ in $C$ by
 $$ x_{n+1} =  \frac{\alpha_n}{t_n} \int_{0}^{t_n} T(s) x_n \; d s
  +(1 - \alpha_n) x_n $$
 for $n \in \NaturalNumber$,
 where $\{ \alpha_n \}$ is a sequence in $[0,1]$ such that
 $0 < \liminf_{n} \alpha_n \leq \limsup_{n} \alpha_n \allowbreak < 1$,
 and $\{ t_n \}$ be a sequence in $(0,\infty)$
 satisfying
 $ \lim_{n} t_n = \infty$ and
 $\lim_{n} t_{n+1} / t_n = 1 $.
Then $\{ x_n \}$ converges strongly to
 a common fixed point $z_0$ of $\{ T(t) : t \geq 0 \}$.
\end{cor}

\begin{cor}[\cite{REF:TSP_retopit_1_03}]
\label{COR:chara-R+}
Let $E$ be a Banach space with the Opial property and
 let $C$ be a weakly compact convex subset of $E$.
Let $\{ T(t) : t \geq 0 \}$ be
 a one-parameter strongly continuous semigroup of nonexpansive mappings on $C$.
Then for $z \in C$,
 the following are equivalent:
\begin{enumerate}
\item
 $z$ is a common fixed point of $\{ T(t) : t \geq 0 \}$;
\item
 there exists a subnet of
 a net
 $$ \left\{ \frac{1}{t} \int_{0}^t T(s)x \; ds \right\} $$
 in $C$
 converging weakly to $z$.
\end{enumerate}
\end{cor}

We also obtain the following new results.

\begin{cor}
\label{COR:chara-mu-R+-Opial}
Let $C$ be a weakly compact convex subset of
 a Banach space $E$ with the Opial property.
Let $\{ T(t) : t \geq 0 \}$ be
 a one-parameter strongly continuous semigroup of nonexpansive mappings on $C$.
Let $X$ be as in Lemma \ref{LEM:R+} and
 let $\mu$ be an invariant mean on $X$.
Then
 $z \in C$ is a common fixed point of $\{ T(t) : t \geq 0 \}$
 if and only if $T_\mu z = z$.
\end{cor}

\begin{cor}
\label{COR:chara-mu-R+-compact}
Let $C$ be a compact convex subset of a Banach space $E$.
Let $\{ T(t) : t \geq 0 \}$ be
 a one-parameter strongly continuous semigroup of nonexpansive mappings on $C$.
Let $X$ be as in Lemma \ref{LEM:R+} and
 let $\mu$ be an invariant mean on $X$.
Then
 $z \in C$ is a common fixed point of $\{ T(t) : t \geq 0 \}$
 if and only if $T_\mu z = z$.
\end{cor}

\section{Appendix}
\label{SC:appendix}

In this section,
 using the notion of universal nets,
 we give an invariant mean.
We recall that a net $\{ y_\alpha : \alpha \in D \}$
 in a topological space $Y$ is universal if
 for each subset $A$ of $Y$,
 there exists $\alpha_0 \in D$ satisfying
 either of the following:
 \begin{itemize}
 \item $y_\alpha \in A$ for all $\alpha \in D$ with $\alpha \geq \alpha_0$; or
 \item $y_\alpha \in Y \setminus A$ for all $\alpha \in D$
 with $\alpha \geq \alpha_0$.
 \end{itemize}
For every net $\{ y_\alpha : \alpha \in D \}$,
 by the Axiom of Choice,
 there exists a universal subnet $\{ y_{\alpha_\beta} : \beta \in D' \}$
 of $\{ y_\alpha : \alpha \in D \}$.
If $f$ is a mapping from $Y$ into a topological space $Z$ and
 $\{ y_\alpha : \alpha \in D \}$ is a universal net in $Y$,
 then $\{ f(y_\alpha) : \alpha \in D \}$ is also a universal net in $Z$.
If $Y$ is compact, then
 a universal net $\{ y_\alpha : \alpha \in D \}$ in $Y$ always converges.
See \cite{REF:Kelley_TOP} and others for details.

\begin{prop}
\label{PROP:invariant-mean}
Let $S$ be a commutative semigroup,
 let $X$ be a linear subspace of $B(S)$
 such that $I_S \in X$ and $X$ is $\ell_s$-invariant for all $s \in S$.
Let $\{ \mu_\alpha : \alpha \in D \}$ be
 an asymptotically invariant net of means on $X$.
Let $\{ \mu_{\alpha_\beta} : \beta \in D' \}$ be a universal subnet of
 $\{ \mu_\alpha : \alpha \in D \}$.
Define a function $\mu$ from $X$ into $\RealNumber$ by
 $$ \mu(a)
 = \lim_{\beta \in D'} \mu_{\alpha_\beta} (a) $$
 for all $a \in X$.
Then $\mu$ is an invariant mean on $X$.
\end{prop}

\begin{proof}
Since the net $\{ \mu_{\alpha_\beta} : \beta \in D' \}$ is universal,
 $\mu$ is well-defined.
Since $\{ \mu_{\alpha_\beta} : \beta \in D' \}$ is
 also an asymptotically invariant net of means on $X$,
 $\mu$ is an invariant mean;
 see \cite{REF:Day1959_Illinois}.
This completes the proof.
\end{proof}

\end{document}